\documentclass[11pt]{article}
\usepackage{amssymb}
\usepackage{graphicx}
\newcommand{\beq}{\begin{equation}}
\newcommand{\bea}{\begin{eqnarray}}
\newcommand{\beas}{\begin{eqnarray*}}
\newcommand{\eeq}{\end{equation}}
\newcommand{\eea}{\end{eqnarray}}
\newcommand{\eeas}{\end{eqnarray*}}
\def\begeq{\begin{equation}}
\def\endeq{\end{equation}}
\newcommand{\benu}{\begin{enumerate}}
\newcommand{\eenu}{\end{enumerate}}
\newcommand{\bit}{\begin{itemize}}
\newcommand{\eit}{\end{itemize}}

\begin{document}
\title{DE$-$SUSPENSION OF FREE $S^3-$ACTIONS
ON HOMOTOPY SPHERES}
\vspace{4cm}
\author{Issam. H. Kaddoura.\footnote{e.mail~:~issam.kaddoura@liu.edu.lb.}~\footnote{
Faculty of Science Baghdad University. } \\
Lebanese International University \\Department of Mathematics
 }
\date{}
\maketitle
\begin{abstract}
In the paper of Montgomery, $D$. and Yang, $C.T.$ \cite{5}, they discuss the de$-$suspension of
smooth free actions of $S^1$ on $(2n+1)-$dimensional homotopy spheres. In this paper we discuss the
de-suspension of smooth free actions of $S^3$ on $(4n+3)-$dimensional homotopy spheres.
\end{abstract}
\section{Introduction }
Through this paper, $R^n$ denotes the Euclidean $n-$space, $S^n$ denotes the unit $n-$sphere in
$R^{n+1}$ and $QP(n)$ the quaterionic projective space, all having the usual differentiable structures.
By a homotopy $n-$sphere abbreviated by $\sum^n$, we mean a closed differentiable $n-$manifold having
the homotopy type of $S^n$, and by a homotopy quaterionic projective $n-$space, abbreviated by
$HQP(n)$, we mean a closed differentiable $4n-$manifold having the homotopy type of $QP(n)$.
$\large{\pi_n} (M)$ denotes the $n^{th}-$homotopy group of $M$, $H_i(M)$ and $H^i(M)$ denote the homology and
cohomology of a space $M$, and assumed to be satisfying the Eeilenberg$-$Steenrod axioms, See
[\cite{13},p.6]. It is well known that $S^1$ and $S^3$ are the only compact connected Lie groups which have
free differentiable actions on homotopy spheres, \cite{8,6}. It follows from Gleason's lemma \cite{1} that
such an action is always a principal fibration which is homotopically equivalent to the classical
Hopf fibration. In fact, there are always infinitely many differentiably distinct free actions of $S^3$
on $\sum^{4n+3}$ for $n \geq 2$, \cite{16}, such that for any differentiable action of the group $S^3$ on$\sum^{4n+3}$, the
orbit space $\sum^{4n+3}/S^3$ is a $HQP(n)$, and that for any $HQP(n)$ there is, up to differentiable
equivalence a unique free differentiable action of $S^3$ on $\sum^{4n+3}$, such that $\sum^{4n+3}/S^3$ is
diffeomorphic to $HQP(n)$.\\
Hence, there is one to one correspondence between differentiable
equivalence classes of free differentiable actions of $S^3$ on $\sum^{4n+3}$ and diffeomorphism classes of
$HQP(n)'s$. If we let $M=\sum^{4n+3}/S^3$, there is a differentiable homotopy equivalence $f :M \rightarrow
QP(n)$ which is transverse regular on $QP(n-1)$, denoted $f\pitchfork QP(n-1)$. Therefore $N=f^{-1}(QP(n-1))$ is a differentiable
 submanifold of $M$ of codimension $4$. By means of framed surgery
as devised in Browder \cite{15} and Novikov \cite{9}, it is always possible to kill the kernel of\\
$f_{\ast i} : \pi_i (N)\rightarrow \pi_i (QP(n-1))~~~\mbox{for}~~i=0,1,\ldots,2n-3.$(below the middle dimension).
That means, there is a differentiable homotopy equivalence $f^1: M\rightarrow QP(n)$ which is homotopic to $f$ and
$f^1 \pitchfork QP(n-1)$ such that if we let $N^1= f^{-1}(QP(n-1))$, \\then ${f}^1_{\ast i}:\pi_i(N^1)\rightarrow \pi_i(QP(n-1))$
has a trivial kernel for $i= 0,1,2,\ldots 2n-3$. But for $i=2n-2$ (at the middle dimension) there is an
obstruction called the index. In fact a necessary and sufficient condition for killing the kernel of
$f_{\ast i}$ for $i=2n-2$, is that this index vanishes. This index turns out to be dependent only on the action
$(S^3,\sum^{4n+3})$ see \cite{7}. We note that if we succeed to kill the kernel of $f_{\ast i}:\pi_i(N)\rightarrow \pi_i(QP(n-1))$  for
$i= 0,1,2,\ldots 2n-2$, hence by Poincare duality on the kernels,we succeed to kill the kernel of $f_{\ast i}$ for all
$i= 0,1,2,\ldots 4n-4$, then there is a differentiably imbedded $HQP(n-1)$ in $HQP(n)$ and hence there is a
differentiably imbedded $S^3-$ invariant $\sum^{4n-1}$ in $\sum^{4n+3}$ called a characteristic homotopy $(4n-1)-$
sphere this what is called a de$-$suspension of the action.
\section{Basic Construction }
In this section, we give a construction which will be needed in our work. It is essentially
the same construction as seen in \cite{15} and \cite{9} with certain modifications are made as in Haefliger
\cite{2}, Levine \cite{9}, Montgomery and Yang \cite{5}. For any manifold $M$ ,$\tau(M)$ denotes the tangent
bundle of $M$ and if $N$ is a submanifold of $M$ , $\nu(N,M)$ denotes the normal bundle of $N$ in $M$. By an
$r-$frame on $N$, we mean an $r-$tuple of linearly independent vector fields in $\tau(M)|_{N}$ . $D^n$ denotes
the closed unit disk in $R^n$ and $B^n= D^n - S^{n-1}$. Let $M$ and $M^1$ be closed $\ell -$ manifolds, $N^1$ a closed
$m-$submanifold of $M^1$, $f : M \rightarrow M^1$ a map which is transverse regular on $N^1$, $N= f^{-1}(N^1)$ and $h :
S^{k-1} \rightarrow N$ is an imbedding. We now state assumptions as by Haefliger in \cite{2} and Levine in \cite{9}
and exactly as Montgomery and Yang in \cite{5}. Later these assumptions will be verified for the
construction made in this paper.
\begin{itemize}
\item[(A)] $ 0 \leq k-1 \leq\frac{m}{2}$  and  $2 \leq m + 2 \leq\ell.$
\item[(B)] $f\circ h: S^{k-1} \rightarrow N^1$ is null homotopic and in fact $f\circ h: (S^{k-1})= y$ a point in $N^1$.
\item[(C)] $\pi_i (M - N) = 0~~\mbox{for}~~ i\leq k-1.$
\item[(D)] $\pi_i (M^1 - N^1) = 0~~\mbox{for}~~ i\leq \ell.$
\end{itemize}
Choose an $(\ell-m)$ frame of $\nu(M^1-N^1)$  at $y$ and let $(v_1,v_2,\ldots,v_{\ell-m})$ be its inuced $(\ell-m)$ frame on
$h(S^{k-1})$ by $f$. Then $h$ can be extended to an imbedding of $D^k$ into $M$, also denoted by $h$, such that
$h(B^k ) \subset M-N$ and $v_1$ points radically into $h(D^k)$. we also assume :
\begin{itemize}
\item[(E)]  The frame $(v_1,v_2,\ldots,v_{\ell-m})$ can be extended to an $(\ell-m-1)$ frame on $h(D^k )$ in $\nu(h(D^k ), M)$.
The obstruction to this extention is an element $\lambda \in \pi_{k-1}(V_{\ell-m-1}(R^{\ell-k}))$  where $V_{\ell-m-1}(R^{\ell-k})$ is the
stiefel manifold of orthonormal $(\ell-m- 1)$ frames in $R^{\ell-k}$  [\cite{4},p.83]. Hence $(E)$ means $\lambda=0$. Notice
that when $h :S^{k-1}\rightarrow V_{\ell-m-1} (R^{\ell-k} )$, then $h$ can be extended to an imbedding $h^1: D^k \rightarrow V_{\ell-m-1}(R^{\ell-k} )$,
provided the obstruction which is an element $\lambda\in H^k (D^k ,S^{k-1}, \pi_{k-1}(V_{\ell-m-1}(R^{\ell-k} )))$ vanishes.
Using the exact cohomology sequence of the pair $(D^k, S^{k-1})$:
\end{itemize}
\begin{itemize}
\item[]$\rightarrow H^{k-1}(D^k ,G)\rightarrow H^{k-1}(S^{k-1},G) \rightarrow H^k (D^k ,S^{k-1},G) \rightarrow H^k (S^{k-1},G) \rightarrow$
\item[]where $G=\pi_{k-1}(V_{\ell-m-1}(R^{\ell-k} ))$ , we get :
\item[]$0 \rightarrow \pi_{k-1}(V_{\ell-m-1}(R^{\ell-k} ))\rightarrow H^k (D^k ,S^{k-1},\pi_{k-1}(V_{\ell-m-1}(R^{\ell-k} ))) \rightarrow 0$
\item[]Hence, $H^k (D^k ,S^{k-1},\pi_{k-1}(V_{\ell-m-1}(R^{\ell-k} )))= \pi_{k-1}(V_{\ell-m-1}(R^{\ell-k} ))$.
\end{itemize}
Moreover, if $2(k-1) \leq m,$ then $\pi_{k-1}(V_{\ell-m-1}(R^{\ell-k} ))\cong 0$,
  i.e, $(E)$ holds for $2(k-1) \leq m$. Under these
assumptions Haefliger in \cite{2} has shown that the map $f$ can be extended to a map $F: M \times[0,1] \rightarrow M^1$
such that $F\pitchfork N^1$ and $B=F^{-1}(N^1)$ is a compact $(m +1)-$manifold in $M\times[0,1]$ which is a
framed cobordism from $N$ to $N^1= B \cap (M\times\{1\})$ obtained by attaching a handle to $N\times[0,1]$ at
$h(S^{k-1})$. Let $f_1= F|_{M\times\{1\} }$, which is also regarded as the map of $M$ into $M^1$ given by $x\rightarrow F(x,1).$
Then $f_1$ is a homotopy equivalence such that $f \sim f_1,f_1 \pitchfork N^1$ and $N_1= f_1^{-1}(N^1)$ is obtained from
$N=f^{-1}(N^1)$ by performing a framed surgery at $h(S^{k-1})$.
\section{ Surgery below the middle dimension }
Suppose that a free differentiable action of $S^3$ on $\sum^{4n+3}$ is given, $n \geq 2$.\\
Let $M=\sum^{4n+3}/S^3$,
then there is a homotopy equivalence $f: M \rightarrow QP(n)$, let $\sigma: S^i \rightarrow M$ be an imbedding
representing a generator of $\pi_i (M)$, by the relative transversality theorem, we may assume $f$ is
transverse regular on $QP(n-1)$, hence $N= f^{-1}(QP(n-1))$ is a closed manifold of $M$ of codimension
$4$ and containing $\sigma(S^i)$ for $i < 2n-1$.\\\\
{\bf  Lemma 3.1 } $~$
\begin{itemize}
\item[(a)] $M- N$ is connected, i.e., $\pi_0(M-N)= 0.$
\item[(b)] $\pi_i(M-N)= 0$ for $i=1, 2.$
\end{itemize}
{\bf Proof.} $~$ \\
$(a)$ let $x,y \in M-N$ then $x,y\in M,$ but $M$ is connected then there exist a path $\alpha: [0,1] \rightarrow M$ such
that $\alpha(0)=x$ and $\alpha(1)=y$. By transversality homotopy theorem, see [\cite{14},p.70], there exist
a path $\alpha^1:[0,1]\rightarrow M$ such that $\alpha^1(0)=x$ and $\alpha^1(1)=y$ and $\alpha^1 \pitchfork N$.
But $dim\alpha^1([0,1]) +dimN < dimM$, thus $\alpha^1 \cap N=\phi $ and $\alpha^1([0,1]) \subset M -N,$ hence $M-N$ is
pathwise connected and therefore is connected.\\
\\
$(b)$ Let $\alpha : S^1\rightarrow M-N$ be an imbedding representing a generator of $\pi_1(M-N)$, $\alpha$ can be
regarded as an imbedding into $M,$ $\alpha: S^1\rightarrow M.$ But $\pi_1(M)=0$ then $\alpha$ can be extended to
$\alpha^1: D^2\rightarrow M$ such that $\alpha^1|_{\partial D^2}= \alpha$, $\alpha$ is homotopic to $\alpha^1$ and $\alpha^1\pitchfork N.$
But dim$ \alpha^1(D^2) + \mbox{dim}N < \mbox{dim}M$ thus $\alpha^1(D^2)\cap N =\phi$, which implies $\alpha^1(D^2)\subset M-N$
therefore $\pi_1(M-N) =0$, similarly $\pi_2(M-N) =0.$\\\\
{\bf Lemma 3.2 }: $~$ $f_{\ast i}: H_i(N)\rightarrow H_i (QP(n-1)) $ is onto $\forall i$.\\
In particular $f_{\ast4n-4}: H_{4n-4}(N)\rightarrow H_{4n-4}(QP(n-1))$ is an isomorphism.\\
{\bf proof.} $~$ Let $A^1$ be the complement of an open tubular neighborhood of $QP(n-1)$ in $QP(n)$ such that
$A=f^{-1}(A^1)$ is the complement of open tubular neighborhood of $N$ in $M$. It is clear that
$f^{\ast}: H^{4n} (QP(n)) \rightarrow H^{4n} (M)$ is an isomorphism. Consider the cohomology exact sequence of the
pair $(QP(n),A^1):$\\
\[
\rightarrow H^{4n-1}(A^1) \rightarrow H^{4n} (QP(n),A^1) \rightarrow H^{4n} (QP(n)) \rightarrow H^{4n} (A^1)\rightarrow
\]
But, $QP(n)=QP(n-1)\cup D^{4n},~ A^1$ is a deformation retract of $QP(n) - QP(n-1)$, thus $A^1$ has the
homotopy type of $D^{4n}$, hence $H^{4n-1}(A^1) = 0$ and $H^{4n} (A^1) = 0.$
Therefore $H^{4n} (QP(n),A^1) \rightarrow H^{4n} (QP(n))$ is an isomorphism. Moreover,the cohomology exact
sequence of $(M,A)$:\\
\[
\rightarrow H^{4n-1}(A) \rightarrow H^{4n} (M,A) \rightarrow H^{4n} (M) \rightarrow H^{4n} (A) \rightarrow
\]
$A$ is a deformation retract of $M-N$, and by Hurewicz theorem, [\cite{10},p.349]:\\
$H_1(A)\cong \pi_1(A)\cong \pi_1(M-N)=0.$ We get $H^{4n-1}(A)=0$ and $H^{4n}(A)= 0.$
Therefore $H^{4n}(M,A) \rightarrow H^{4n} (M)$ is an isomorphism. By the isomorphism theorem of Thom, see
\cite{6}, $H^{4n-4}(QP(n-1)) \rightarrow H^{4n} (QP(n),A^1)$ and $H^{4n-4}(N) \rightarrow H^{4n} (M,A)$ are isomorphisms. \\It follows
that $f^{\ast}: H^{4n-4}(QP(n-1)) \rightarrow H^{4n-4}(N)$ is an isomorphism.\\\\
Claim: $f^{\ast i}: H^{i}(QP(n-1))  \rightarrow H^i (N)$ is one to one $\forall~ i$.\\
{\bf Proof.} $~$ Let $\alpha\neq 0 \in H^i(QP(n-1))$, then	$ \exists~~ \beta \in H^{4n-4-i} (QP(n-1))$ such that \\
$\alpha \cup \beta=d\neq 0\in  H^{4n-4}(QP(n-1))$. This implies $f_{4n-4}^{\ast}(\alpha \cup \beta)=f^{\ast}(\alpha)\cup f^{\ast}(\beta)=f^{\ast}(d)\neq0.$
,i.e., $f^{\ast i}: H^{i}(QP(n-1))  \rightarrow H^i (N)$. is one to one $\forall~ i$. Then by duality between homology and
cohomology $f_{\ast i}:H_i (N)\rightarrow H_{i}(QP(n-1))$ is onto $\forall~ i$\\\\
{\bf Lemma 3.3} $~$ If $k \leq 2n-1$ such that  $f_{\ast i}: H_i (N)\rightarrow H_{i}(QP(n-1)) $ is $1-1$ for $i < k -1$, then $H_{k-1}(M-N) = 0$.\\
{\bf Proof.} $~$ By lemma 2.2 and by the hypothesis we get  $f_{\ast i}: H_i (N)\rightarrow H_{i}(QP(n-1)) $ is an isomorphism
for $i < k-1$. Therefore, by poincare duality [\cite{12},p.138],\\
 $f^{\ast }: H^{4n-4-i}(QP(n-1)) \rightarrow  H^{4n-4-i} (N)$ is an
isomorphism for $i < k-1$, so that $H^{4n-4-i} (M) \rightarrow H^{4n-4-i}(N)$ induced by the inclusion map is an
isomorphism for $i < k-1$, \cite{12}. Consider the exact cohomology sequence of $(M,N)$:\\
\[ \rightarrow H^{4n-4-i}(M,N) \rightarrow H^{4n-4-i}(M) \rightarrow H^{4n-4-i} (N)\rightarrow H^{4n-3-i} (M,N) \rightarrow\]
$\mbox{for}~~ i = k-4:$
\[\rightarrow H^{4n-k}(M,N) \rightarrow H^{4n-k}(M) \rightarrow H^{4n-k} (N)\rightarrow H^{4n-k+1} (M,N) \rightarrow \]
$H^{4n-4}(M) \cong H^{4n-k}(N)$ and  $H^{4n-k+1} (M) \cong H^{4n-k+1} (N) $ hence, $ H^{4n-k+1} (M,N)\cong 0$\\
by Lefshetz duality \cite{11}, we get $H_{k-1}(M-N)\cong H^{4n-k+1}(M,N)\cong 0$.
\\\\\\
{\bf Lemma 3.4} $~$ If $k \leq 2n-1$ such that $f_{\ast i}:\pi_i (N)\rightarrow \pi_{i}(QP(n-1))$ is $1-1$ for $i < k-1$, then both
$f_{\ast i}:\pi_i (N)\rightarrow \pi_{i}(QP(n-1))$, and $f_{\ast i}:H_i (N)\rightarrow H_{i}(QP(n-1))$ are isomorphisms for $i < k-1$.
Moreover, if $\rho:\pi_{k-1}(N)\rightarrow H_{k-1}(N)$ is the Hurewicz homomorphism, $ker(f_{\ast \pi_{k-1}} )$ is the kernel of
$f_{\ast }:\pi_{k-1} (N)\rightarrow \pi_{k-1}(QP(n-1))$ and  $ker(f_{\ast H_{k-1}} )$ is the kernel of $f_{\ast }:H_{k-1} (N)\rightarrow H_{k-1}(QP(n-1))$, then for $k \neq 2~~ \rho(ker(f_{\ast \pi_{k-1}} ))\cong ker(f_{\ast H_{k-1}} )$ and for $k= 2,~~\rho$ is onto and its kernel is the commutator
subgroup of $ker(f_{\ast \pi_{k-1}} )$. Furthermore $\pi_{i} (M- N) = 0$ for $i\leq k-1.$\\
{\bf Proof. }  $~$ Let $\sigma: S^i \rightarrow M$ be an imbedding representing a generator of $\pi_i (M),~\sigma(S^i)$ is contained
in $N$ as our assumption at the beginning of this section for $i < 2n-1,$  $f{\rho}: S^i \rightarrow QP(n-1)$
represents a generator of $\pi_i (QP(n-1))$ for $i\leq k-1$, so that $f_{\ast i}:\pi_i (N)\rightarrow \pi_{i}(QP(n-1))$ is an
isomorphism for $i< k-1.$ and onto for $i\leq k-1.$. Let $C_f$ be the mapping cylinder of $f : N \rightarrow QP(n-1)$
 ,i.e., $C_f= N\times I \cup_{f} QP(n-1)$, it is clear that $QP(n-1)$ is homotopically equivalent to $C_f$. See[\cite{11},p.283].
We consider the two exact sequences which together with the Hurewicz homomorphism form the
commutative diagram:
\begin{eqnarray*}
\rightarrow \pi_i (N) \buildrel {f_{\ast }} \over\longrightarrow  ~&\pi_{i}(QP(n-1))&~\rightarrow \pi_{i}(C_f,N)\rightarrow \pi_{i-1} (N)\rightarrow \hspace{3cm}\\
\rho\downarrow~~~~~~&\rho\downarrow&~~~~~~~\rho\downarrow~~~~~~~~~~\rho\downarrow\hspace{3cm}\\
\rightarrow H_i (N) \longrightarrow~&H_{i}(QP(n-1))&~\rightarrow H_{i}(C_f,N)\rightarrow H_{i-1} (N)\rightarrow \hspace{3cm}
\end{eqnarray*}
It is clear that $\pi_i (C_f,N)= 0$ for $i \leq k-1$, hence $H_i(C_f ,N)=0$
 for $i \leq k-1$, and $\pi_k (C_f,N)\cong H_k (C_f,N)$ ( By Hurewicz theorem ),
 \cite{3}. Moreover,\\
 $f_{\ast i}:H_{i} (N)\rightarrow H_{i}(QP(n-1))$ is an isomorphism for $i< k-1,$
and onto for $i\leq k-1$.\\
\\
Now, the following diagram is commutative and its rows are exact :
\begin{eqnarray*}
~\rightarrow \pi_{k}(C_f,N)\rightarrow ~&ker(f_{\ast \pi_{k-1}})~&\rightarrow 0 \hspace{3cm}\\
~~~~~~~\rho\downarrow~~~~~~~~~~&\rho\downarrow~&~~~~~~~~~~\hspace{3cm}\\
~\rightarrow H_{k}(C_f,N)\rightarrow ~&ker(f_{\ast H_{k-1}} )~& \rightarrow 0\hspace{3cm}
\end{eqnarray*}
We get $\rho (ker(f_{\ast \pi_{k-1}}))\cong ker(f_{\ast H_{k-1}})$, see\cite{10}.
since $f_{\ast i}:H_{i} (N)\rightarrow H_{i}(QP(n-1))$ is an isomorphism for $i< k-1.$
and onto for $i\leq k-1$, \cite{13}, so as $H_i (N) \rightarrow H_i (M)$ induced by the inclusion map,
similarly as in (Lemma 3.3) we get $H_i (M-N) = 0 $ for $i\leq k-1.$ But $\pi_0 (M-N) = 0 $
and $\pi_1 (M-N) = 0 $ (By Lemma 3.1), hence by Hurewicz theorem $\pi_i (M-N) = 0 $ for $i\leq k-1.$\\\\
Now, we are in position to use the construction of the preceding section to kill the kernel of
$f_{\ast i}:\pi_{i} (N)\rightarrow \pi_{i}(QP(n-1))$ for $i<2n-2$.
\\Consider $f_{\ast }:\pi_{0} (N)\rightarrow \pi_{0}(QP(n-1))$, if the kernel is not
trivial  ,i.e., $N$ has more than one component, so there is an imbedding $h : S^0 \rightarrow N$ such that $h(S^0)$ is
not contained in the same component of $N$. If we put $M,~ QP(n),~ QP(n-1),~ f~\mbox{and}~ h$ in place of $ M,~
M^1,~ N^1,~ f ~\mbox{and}~ h$ respectively. We have to satisfy the required assumptions,then $(A),(B),(C),(D)$ and $(E)$ are easily
verified. Therefore, we obtain a homology equivalence $f_1: M \rightarrow QP(n)$ such that $f \sim f_1$ and
$f_1 \pitchfork QP(n-1)$ such that $f_1^{-1}(QP(n-1))$ has one less component than $f^{-1}(QP(n-1))$.
Repeating this construction, if necessary, we can finally have a homotopy equivalence $f^1: M \rightarrow QP(n)$ such that
$f \sim f^1,~ f^1 \pitchfork QP(n-1)$ and ${f^1}^{-1} (QP(n-1))$ is connected and thus the kernel of $f_{\ast }^1:\pi_{0}({f^1}^{-1}(QP(n-1)) \rightarrow \pi_0 (QP(n-1)) $  is trivial, that is we kill the kernel of
$f_{\ast }:\pi_{0}(N) \rightarrow \pi_0 (QP(n-1)) $ by using $f^1$ to replace $f$. Consider $f_{\ast }:\pi_{1}(N) \rightarrow \pi_1 (QP(n-1)) $ such
that the kernel of $f_{\ast }:\pi_{0}(N) \rightarrow \pi_0 (QP(n-1)) $ is already killed.
Given an element $\mu$ of the kernel of $f_{\ast }:\pi_{1}(N) \rightarrow \pi_1 (QP(n-1)) $, there is an imbedding $h : S^1 \rightarrow
N$ representing $\mu$ with $M,~ QP(n),~QP(n-1),~f ~\mbox{and}~ h$ in place of $M,~ M^1,~N^1,~ f ~\mbox{and}~ h$ respectively. we
have to satisfy the required assumptions: $(A),(B),(C),(D)$ and $(E)$ for $k=2,$ $m = 4n-4 $ and $\ell= 4n$.
\begin{itemize}
\item[(A)] It is obvious that $0 \leq k-1 \leq \frac{m}{2}$ and $5 \leq m+2 \leq \ell$.
\item[(B)] $\mu \in kerf_{\ast}$ implies $f_{\ast}(\mu) = 0$, hence $f\circ h: S^1 \rightarrow QP(n-1)$ is null homotopic.
\item[(C)] $\pi_0 (M-N) = 0 $ and $\pi_1 (M-N) = 0 $ (by lemma 3.1).
\item[(D)] $\pi_i (QP(n)-QP(n-1)) = 0,~~\forall~~ i \leq 4n .$
\end{itemize}
In fact, $QP(n)- QP(n-1)$ has the homotopy type of $D^{4n}$. Moreover $h : S^1 \rightarrow
N$ can be extended to an imbedding $h: D^2 \rightarrow M$ such that $h(B^2) \subset M-N$
and $v_1$ points radically into $h(D^2)$, where $v_1$ as in section $2$.
\begin{itemize}
\item[(E)] $\mu \in 	 \pi_1 (V_3(R^{4n-2})) = 0.$ See [\cite{9},p.83].
\end{itemize}
Hence, there is a homotopy equivalence $f_1: M \rightarrow QP(n)$ such that $f \sim f_1,~ f_1 \pitchfork QP(n-1)$ and such
that the kernel of $f_{1 \ast }:\pi_{0}(f_1^{-1}(QP(n-1))) \rightarrow \pi_0 (QP(n-1)) $ is trivial and the kernel of
$f_{1 \ast }:\pi_{1}(f_1^{-1}(QP(n-1))) \rightarrow \pi_1 (QP(n-1)) $ is the quotient group of the kernel of
$f_{ \ast }:\pi_{1}(f^{-1}(QP(n-1))) \rightarrow \pi_1 (QP(n-1)) $ by the normal subgroup generated by $\mu$.
Repeating the same procedure if necessary, the kernel of $f_{ \ast }:\pi_{1}(f^{-1}(QP(n-1))) \rightarrow \pi_1 (QP(n-1)) $
can be killed also. Similarly, if $k$ is an integer such that $1 < k-1 < 2n-2$ and the kernel of
$f_{ \ast }:\pi_{i}(N) \rightarrow \pi_i (QP(n-1)) $ has been killed for $i < k-1$, we can kill the kernel of
$f_{ \ast }:\pi_{k-1}(N) \rightarrow \pi_{k-1} (QP(n-1)) $ by using the same argument. Finally, we have a homotopy
equivalence $f^1: M \rightarrow QP(n)$ such that $f \sim f^1,~ f^1 \pitchfork QP(n-1)$ and the kernel of
$f_{\ast i}^1:\pi_{i}({f^1}^{-1}(QP(n-1))) \rightarrow \pi_i (QP(n-1)) $ is trivial for $i < 2n-2.$
\section{The  obstruction at the middle dimension (The index)
}
Assume that a free differentiable action of the group $S^3$ on $\sum^{4n-3}$ is given, $n \geq 2$ and let
$M= \sum^{4n-3}/S^3$ and $f : M \rightarrow QP(n)$ is a homotopy equivalence which is transverse regular on
$QP(n-1)$ and let $N = f^{-1}(QP(n-1))$. As in the preceding section, we may assume that the kernel
of $f_{ \ast i}:\pi_{i}(N) \rightarrow \pi_i (QP(n-1)) $ has been killed for $i < 2n-2$.
Therefore by lemma 3.4,  $f_{ \ast i}:\pi_{i}(N) \rightarrow \pi_i (QP(n-1)) $ and  $f_{ \ast i}:H_{i}(N) \rightarrow H_i (QP(n-1)) $ are
isomorphisms for $i < 2n-2$ and the Hurrewicz homomorphism $\sigma$
 maps $ker(f_{\ast \pi_{2n-2}} )$
isomorphically onto $ker(f_{\ast H_{2n-2}} )$. Consider the bilinear form $F : ker(f_{\ast H_{2n-2}} )\otimes ker(f_{\ast H_{2n-2}} )\rightarrow Z$
which maps each $\xi \otimes \mu \in ker(f_{\ast H_{2n-2}} )\otimes ker(f_{\ast H_{2n-2}} )$ into the intersection number $\xi.\mu$. It is symetric
and has a signature $I$. The integer $I$ is independent of the choice of $f$. For more details see \cite{5}.
$I$ is called the index of the action $(S^3,\sum^{4n+3})$ and is written as $I(S^3,\sum^{4n+3})$.
\\
\\
{\bf Lemma 4.1}  $~$ Suppose $n= 2s$ and the group $ker(f_{\ast \pi_{s}} )$ is trivial for $i < s$. Then $ker(f_{\ast H_{s}} )$ is free abelian, and the
matrix of intersections of the base cycles of the group $ker(f_{\ast H_{s}} )$ is unimodular. [\cite{11},p.283].\\\\
{\bf Theorem 4.1} $~$ For any integer $n \geq 2$, the index $I$ vanishes iff there is a homotopy equivalence $f : M \rightarrow QP(n)$
which is transverse regular on $QP(n-1)$ and such that if $N = f^{-1}(QP(n-1)),$ then $f : N \rightarrow QP(n-1)$
is a homotopy equivalence.\\
{\bf Proof.} $~$ If $f : M \rightarrow QP(n)$ is a homotopy equivalence such that $f :N \rightarrow QP(n-1)$ is a homotopy
equivalence then $ker(f_{\ast H_{i}} )=0~~\forall~i$ and hence $I$ vanishes. Conversely, assume that $I$ vanishes and
let $f : M \rightarrow QP(n)$ is a homotopy equivalence described at the beginning of this section.
Using lemma 4.1, $ker(f_{\ast H_{2n-2}} )$ is a free abelian group which has the basis
$\{\xi_1,\xi_2,\ldots\ldots,\xi_{2r-1},\xi_{2r}\}$ such that
\[
\xi_i.\xi_j=
\left\{\begin{array}{ll}
1~~~~~\mbox{if}~~ (i,j)=(2t-1,2t)~~ \mbox{for some t}
\vspace{.2cm}\\
0~~~~~~~~~~~~~otherwise
\end{array}
\right.
\]
Let $h : S^{2n-2} \rightarrow N$ be an imbedding representing $\xi_1$. We can use the basic construction of section 2
to kill $\xi_1$ and $\xi_2$ at the same time. For $\ell= 4n,~ k= 2n-1,~ m=4n-4,~\\ M=\sum^{4n+3}/S^3,~ M^1= QP(n),~
N^1= QP(n-1)~~\mbox{and}~~ N= f^{-1}(QP(n-1))$.
\begin{itemize}
\item[(A)] It is obvious that $0\leq  k-1 \leq \frac{m}{2}$ and $5 \leq m +2 \leq\ell .$
\item[(B)] $\xi_1 \in ker(f_{\ast H_{2n-2}} )$, hence $f\circ h : S^{2n-2} \rightarrow N^1$ is null homotopic.
\item[(C)] $\pi_i (M-N) = 0$ for $i \leq 2n-2,$ by lemma 3.4
\item[(D)] $\pi_i (QP(n)-QP(n-1))=0$ for $i \leq 4n.$
\item[(E)] $\lambda \in \pi_{2n-2}(V_3(R^{2n+1}))=0$, [\cite{4},p83].
\end{itemize}
This shows that $\xi_1$ and $\xi_2$ can be killed by the basic construction.\\
Similarly $\xi_3, \xi_4,\ldots\ldots,\xi_{2r-1},\xi_{2r}$
can be killed two at the same time. Hence $ker(f_{\ast H_i} )$ is trivial for $i = 1,2,\ldots 2n-2$, and then $ker(f_{\ast \pi_i} )$
is trivial for $i = 1,2,\ldots 2n-2$. By Poincare duality on the kernel's, we deduce that $ker(f_{\ast \pi_i} )$ is trivial
for $i = 1,2,\ldots,4n-4.$ Therefore, $f : N \rightarrow QP(n-1)$ is a homotopy equivalence.
Finally, if $I(S^3,\sum^{4n+3}) = 0$, then there is a differentiably imbedded $HQP(n-1)$ in $HQP(n)$ and
hence there is a differentiably imbedded $S^3-$invariant $\sum^{4n-1}$ in  $\sum^{4n+3}$ which is the de$-$suspension
of the action.\\\\
{\Large\textbf{Acknowledgement}} $~$ I am very grateful to  professor Adil G. Naoum who guided me in different topics of this paper.


\newpage


\begin{thebibliography}{99}
\bibitem{1} A. Gleason , {\em Spaces with a compact Lie group of transformations},
proc. Amer. Math. Soc. I (1950), 35-43.
\bibitem{2} A. Haefliger, {\em Knotted $(4k-1)$ spheres in $6k$ space,}
Ann. of Math.75 (1962), 452-466.
\bibitem{3} C.R. Maunder, {\em Algebraic topology}. Van Nostrand Reinhold Comp.1970.
\bibitem{4} D. Husemoller, {\em Fiber bundles,}
McGraw-Hill, Book Company, New York,(1966).
\bibitem{5} D.Montgomery, and C.T. Yang, {\em Free differentiable actions
on homotopy spheres.} Proc. Con. on transformation groups,
Springer, New York, (1968) 175-192.
\bibitem{6} H.T. Ku, {\em A note on Semi-free actions on homotopy spheres},
Proc. Amer. Math. Soc.22(1969),614-617.
\bibitem{7} H.T. Ku,and M.C.Ku,{\em  Characteristic spheres of free
differentiable actions of $S^1$ and $S^3$ on homotopy spheres.}
Trans. Amer. Math. Soc.156 (1971), 493-504.
\bibitem{8} H.T. Ku,and M.C.Ku, {\em Free differentiable actions of $S^1$ and $S^3$
on homotopy spheres}. Proc. Amer. Math. Soc.25(1970),864-869.
\bibitem{9} J.Levine, {\em A classification of differentiable knots},
Ann. of Math. ,82 (1965), 15-50.
\bibitem{10} P.J. Hilton, and S.Wylie, {\em Homology Theory,}
Cambridge Univ. Press, London and New York,(1960).
\bibitem{11} S.P. Novikov, {\em Homotopically equivalent smooth manifolds,}
Translations Amer. Math. Soc.48 (1965), 271-396.
\bibitem{12} S.T.Hu, {\em Cohomology theory}, Markham Publishing Comp. (1968).
\bibitem{13} S.T. Hu, {\em Homology theory}, Holden-day, Inc.(1970).
\bibitem{14} V. Guillemin, and A. Pollack, {\em Differential topology,
Prentice-Hall,} Inc, New Jersy, (1974).
\bibitem{15} W.Browder , {\em Homotopy type of differentiable manifolds,
Colloquium on algebraic topology,} Aarhus,(1962),42-46.
\bibitem{16} W.C.Hsiang, {\em A note on free differentiable actions of $S^1$ and $S^3$
on homotopy spheres,} Ann. of Math.83 (1966), 266-272.
\end{thebibliography}
\end{document}